\documentclass[12pt]{article}
\usepackage{graphicx}
\usepackage{amssymb}
\usepackage{amsmath}
\usepackage{epstopdf} 

\newtheorem{Pa}{Paper}[section]
\newtheorem{Tm}[Pa]{{\bf Theorem}}
\newtheorem{La}[Pa]{{\bf Lemma}}
\newtheorem{Cy}[Pa]{{\bf Corollary}}
\newtheorem{Rk}[Pa]{{\bf Remark}}

\newtheorem{Ee}[Pa]{{\bf Example}}
\newtheorem{Dn}[Pa]{{\bf Definition}}
\newtheorem{Pn}[Pa]{{\bf Proposition}}

\newcommand{\eh}{\hfill}\newlength{\sperr}

\newenvironment{proof}{{\settowidth{\sperr}{\bf\rm
Proof}%
\par\addvspace{0.3cm}\noindent\parbox[t]{1.3\sperr}
{\textit{ P\eh r\eh o\eh o\eh f\eh .}}%
}}{\nopagebreak\mbox{}\hfill
$\Box$\par\addvspace{0.3cm}}

\numberwithin{equation}{section}

\newcommand{\E}{\mathrm{e}}
\renewcommand{\d}{\,\mathrm{d}}
\newcommand{\I}{\mathrm{i}}

\def\BC{{\mathbb C}}

\title{Integral triangular operators \\ and Friedrichs model}
\author{Lev Sakhnovich}
\date{}

\begin{document}

\maketitle

\thanks{
99 Cove ave., Milford, CT, 06461, USA. 

 E-mail: lsakhnovich@gmail.com}

 \begin{abstract}
In the present paper we investigate a semi-group
 of  triangular integral operators $V_{\beta}$, which is an analogue
 of the semi-group of the fractional integral operators $J^{\beta}$. With the help of these
 semi-groups, we construct and study two classes
 of triangular Friedrichs models $A_{\beta}$ and $B_{\beta}$, respectively. Using generalized wave 
 operators we prove that  $A_{\beta}$ and $B_{\beta}$ are linearly similar to a self-adjoint operator with absolutely continuous spectrum.
 \end{abstract}



\section{Introduction}

\paragraph{1.} In the present paper we investigate
the triangular integral operators
\begin{equation}
V_{\beta}f=
\frac{1}{\Gamma(\beta)}\int_{0}^{x}E_{\beta}(x-t)f(t)\d t,\quad \Re{\beta}>0,\quad f(x){\in}L^{2}(0,\omega),
\label{1.1}
\end{equation}
where   $\Gamma(z)$ is Euler gamma function,
\begin{equation}
E_{\beta}(x)=\int_{0}^{\infty}\frac{1}{\Gamma(s)}\E^{-Cs}s^{\beta-1}x^{s-1}\d s, \quad C\in \BC,
\label{1.2}\end{equation}
and $\BC$ stands for the complex plane. 
We shall show that the operators $V_{\beta}$ form a semi-group, which is an analogue
 of the semi-group of the fractional integral operators
\begin{equation}J^{\beta}f=\frac{1}{\Gamma(\beta)}\int_{0}^{x}(x-t)^{\beta-1}f(t)\d t,\quad \Re{\beta}>0,\quad f(x){\in}L^{2}(0,\omega).\label{1.3}\end{equation}
 With the help
 of the semi-groups of operators $J^{\beta}$ and $V_{\beta}$ we construct and study two classes
 of triangular Friedrichs models $A_{\beta}$ and $B_{\beta}$, respectively. Using generalized wave operators we prove that  $A_{\beta}$ and $B_{\beta}$ are linearly similar to 
 a self-adjoint operator with absolutely continuous spectrum.

\begin{La}\label{Lemma 1.1} The function $E_{\beta}(x)$ is continuous in the domain $(0,\omega]$ and
\begin{equation}
E_{\beta}(x)=\frac{\Gamma(\beta+1)}{x|\ln (x)|^{\beta+1}}\left( 1+O\left(\frac{1}{\ln (x)}\right) \right) ,\quad
x{\to}0.\label{1.4}
\end{equation}
\end{La}
\begin{proof} 
It follows from \eqref{1.2} that the function $E_{\beta}(x)$ is continuous in the domain $(0,\omega]$ and
\begin{equation}  E_{\beta}(x)=\int_{0}^{1}\frac{1}{\Gamma(s)}\E^{-Cs}s^{\beta-1}x^{s-1}\d s+O(1).\label{1.5}
\end{equation}
Using relations
\begin{equation}\frac{1}{\Gamma(s)}=s+O(s^{2}),\quad \E^{-Cs}=1+O(s) \quad \quad (s{\to}0),
\label{1.6}\end{equation} we obtain
\begin{equation}  
E_{\beta}(x)=\frac{1}{x}\int_{0}^{1}s^{\beta}\big( 1+O(s)\big)\E^{-s|\ln (x)|}\d s+O(1).
\label{1.7}\end{equation}
Equality \eqref{1.7} implies the equality \eqref{1.4}. The lemma is proved.
\end{proof}

\begin{Pn}\label{Proposition 1.2}
The operators $V_{\beta},$ defined by formula \eqref{1.1}, are bounded in  all spaces $L^{p}(0,\omega)$, $p{\geq}1$.
\end{Pn}
\begin{proof}
Indeed, according to \eqref{1.4} we have (see \cite[p. 24]{Sakh4})
\begin{equation}
\big\| V_{\beta}\big\|_{p}{\leq}m(\beta)=
\int_{0}^{\omega}\big| E_{\beta}(x)/{\Gamma(\beta)}\big| \d x<\infty.
\label{1.8}
\end{equation}
This proves the proposition.
\end{proof}
The operators $V_{\beta}$ have a following important property \cite{Sakh5}.
\begin{Tm}\label{Theorem 1.3}
The operators $V_{\beta}$\, $(\Re{\beta}>0)$ defined by formula \eqref{1.1},
form a semi-group, that is, 
\begin{equation}V_{\alpha}V_{\beta}=V_{\alpha+\beta},\quad \Re{\alpha}>0,\quad \Re{\beta}>0.\label{1.9}
\end{equation}
\end{Tm}

\paragraph{2.}
 We introduce the following integro-differential operator
\begin{equation}Rf=-\int_{0}^{x}f^{\prime}(t)\ln (x-t)\d t,\quad f(x){\in}L(0,\omega).
\label{1.10}
\end{equation}
In the book \cite[p. 73]{Sakh4} we proved that the operator $V_1$ is a right inverse of  $R$.
\begin{Pn}\label{Proposition 1.4}
If $\beta=1$ and $C=-\Gamma^{\prime}(1)$  then
\begin{equation}
RV_{1}\varphi =\varphi , \quad \varphi^{\prime}(x){\in}L(0,\omega).
\label{1.11}
\end{equation}
\end{Pn}

\paragraph{3.}
 By $H_{\alpha}$ we denote the space of all function such that $f(x){\in}L^{2}(0,\omega)$
and $f(x)=0$ when $x{\in}[0,\alpha]$. It is easy to see that the spaces $H_{\alpha}$,
\,$0<\alpha{\leq}\omega$, are invariant subspaces of the operator $V_{\beta}$.

\begin{Dn}\label{Definition 1.5}
A bounded  operator $T$ is unicellular if its lattice of invariant subspaces is totally ordered
by inclusion.
\end{Dn}
\begin{Tm}\label{Theorem 1.6}
The operator $V_{\beta}$,\, $0<\beta{\leq}1$, defined by formula \eqref{1.1} is  unicellular.
\end{Tm}
\begin{proof}
Let us consider the inner product
\begin{equation}
\big( V_{\beta}^{n}f,g\big) =\int_{0}^{\infty}\frac{1}{\Gamma(s)}\E^{-Cs}s^{n\beta-1}
\left(\int_{0}^{\omega}\int_{0}^{x}(x-t)^{s-1}f(t)\d t\,\overline{g(x)}\d x\right)\d s .
\label{1.12}
\end{equation}
Assuming that $\big( V_{\beta}^{n}f,g\big) =0$ we
use the set of functions $\psi_n(x)=\E^{-x^{2}/2}H_{2n}(x)$, where $H_{2n}(x)$ are Hermite polynomials. The set of functions  $\psi_n(x)$ is complete in the
 Hilbert space $L^{2}(0,\infty)$. We note that $H_{2n}(x)$ are even polynomials of the degree $2n$. Hence, the set of functions  $\varphi_n(x)=\E^{-x^{2}/2}x^{2n}$ is complete in the Hilbert space $L^{2}(0,\infty)$. Changing the variable $s^{\beta/2}=u$ and taking into account the formulas \eqref{1.12} and
$\big( V_{\beta}^{n}f,g\big) =0$, we have
\begin{equation}
\int_{0}^{\omega}\int_{0}^{x}(x-t)^{s-1}f(t)\d t\,\overline{g(x)}\d x=0,
\quad 0<\beta{\leq}1 .
\label{1.14}
\end{equation}
Relation \eqref{1.14} can be written in the form
\begin{equation}
\int_{0}^{\omega}u^{s-1}\int_{u}^{\omega}f(x-u)\overline{g(x)}\d x\d u=0.
\label{1.15}\end{equation}
It follows from \eqref{1.15} that
\begin{equation}\int_{u}^{\omega}f(x-u)\overline{g(x)}\d x=0.
\label{1.16}\end{equation}
By changing the variables $x=\omega-x_{1},\,u=\omega-u_{1}$ we obtain
\begin{equation}\int_{0}^{u_1}f(u_{1}-x_{1})\overline{g(\omega-x_{1})}\d x_{1}=0.
\label{1.17}\end{equation}
Using well-known Titchmarsh's theorem (\cite[Theorem 152]{Tit}) we receive the following assertion: \, 
 if $f(x){\in}H_{\alpha}$ and $f(x){\notin}H_{\gamma}$, when $\gamma>\alpha$, then
the system $V^{n}_{\beta}f$ is complete in the space $H_{\alpha}$.
This proves the theorem.
\end{proof}
\begin{Rk}Using  Titchmarsh's theorem (\cite[Theorem 152]{Tit}) it is easy to prove
that operator $J^{\beta}$,\, $\beta>0$, is unicellular.
\end{Rk}

\section{Triangular integral operators with logarithmic type kernels}
 \label{sec2}
\paragraph{1.} 
In the present section we shall consider the operators
\begin{equation}
S_{\beta}f={\beta}\int_{0}^{x}f(t)\big|\ln (x-t)\big|^{-\beta-1}\frac{\d t}{x-t},
\quad \Re{\beta}>0,
\label{2.1}\end{equation}
in the space $L^{2}(0,\omega)$,\, $0<\omega<1.$
The operator $S_{\beta}$ is the main term in the expression which defines  the operator  $V_{\beta}$ (see \eqref{1.1}, \eqref{1.2}). 
Thus, by  investigating the operator 
$S_{\beta}$ we also obtain
some properties  of the operator  $V_{\beta}$. 
The operator $S_{\beta}$  is of independent  interest  as well.
Rewrite the operator $S_{\beta}$ in the following form
\begin{equation}
S_{\beta}f=\frac{\d}{\d x}\int_{0}^{x}f(t)\big|\ln (x-t)\big|^{-\beta}\d t,
\quad \Re{\beta}{\geq}0.
\label{2.2}
\end{equation}
We note that the operator $S_{\beta}$ is well-defined by formula \eqref{2.2} in
the case where $\Re{\beta}{\geq}0$.
Let us introduce the functions (see \cite{Sakh3} and \cite{Sakh5})
\begin{align}\label{2.3}
\widetilde{s}(\lambda) & =\int_{0}^{\omega}\E^{\I t\lambda}\big|\ln (t)\big|^{-\beta}\d t,
\quad \Re{\beta}{\geq}0,  \\ 
\label{2.4}
\widetilde{s}_{1}(\lambda) & = \int_{0}^{\omega}\E^{\I t\lambda}\big|\ln (t)\big|^{-\beta}
\big( 1-t/\omega\big)\d t, 
\quad \Re{\beta}{\geq}0.
\end{align}
\begin{Tm}
\label{Theorem 2.1} 
If $\Re{\beta}{\geq}0$ and $0<\omega<1$, then the following asymptotic relations are valid:
\begin{align}\label{2.5}
-\I\lambda\widetilde{s}(\lambda) &=\big(\ln |\lambda|\big)^{-\beta}\big( 1+o(1)\big) -
\E^{\I \lambda\omega}\big( -\ln {\omega}\big)^{-\beta},\quad\lambda{\to}\pm\infty, 
\\ \label{2.6}
-\I\lambda\widetilde{s_1}(\lambda) &=\big(\ln |\lambda|\big)^{-\beta}\big( 1+o(1)\big)
,\quad\lambda{\to}\pm\infty.
\end{align}
\end{Tm}
\begin{proof}
According to the Cauchy theorem we have
\begin{equation}
\int_{\gamma}\E^{\I t\lambda}\big(-\ln (t)\big)^{-\beta}\d t=0,
\label{2.7}\end{equation}where the curve $\gamma$ is depicted by Fig.1.

\begin{figure}[h]
    \centering
    \includegraphics[width=0.4\textwidth]{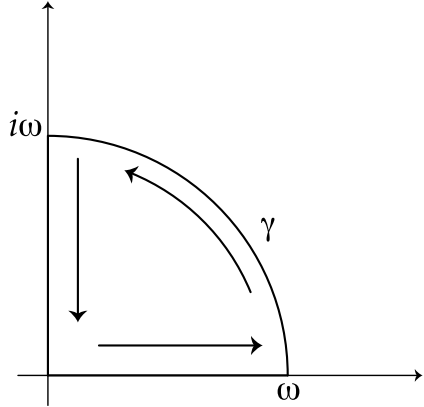}
    \caption{Contour of integration}
    \label{fig:Contour1}
\end{figure}

 From \eqref{2.3}
and \eqref{2.7} we obtain the equality
\begin{equation}
\widetilde{s}(\lambda)=\int_{0}^{\I\omega}\E^{\I t\lambda}\big( -\ln {t}\big)^{-\beta}\d t-
\int_{0}^{\pi/2}\E^{\I t\lambda}\big( -\ln (t)\big)^{-\beta}\I t\d\varphi ,\quad t=\omega \E^{\I \varphi }.
\label{2.8}
\end{equation}
Let us consider the integral
\begin{equation}
\int_{0}^{\I \omega}\E^{\I t\lambda}\big( -\ln {t}\big)^{-\beta}\d t=\int_{0}^{\lambda\omega}
\E^{-u}\big( -\ln {u}+\ln {\lambda}-\I \pi/2\big)^{-\beta}\I \d u/\lambda,\quad\lambda >0.
\label{2.9}
\end{equation}
When $\lambda$ tends to infinity, equality \eqref{2.9} yields
\begin{equation}
\int_{0}^{\I \omega}\E^{\I t\lambda}\big( -\ln {t}\big)^{-\beta}\d t=\frac{\I}{\lambda}\big(\ln \lambda\big)^{-\beta}\big( 1+o(1)\big) ,
\quad \lambda{\to}+\infty.
\label{2.10}\end{equation}
Integrating by parts we estimate the integral
\begin{equation}
\int_{0}^{\pi/2}\E^{\I t\lambda}\big( -\ln (t)\big)^{-\beta}\I t\d\varphi
=-\frac{1}{\I \lambda}\Big( \E^{\I \lambda\omega}\big( -\ln {\omega}\big)^{-\beta}+o(1)\Big) ,
\,\lambda{\to}+\infty,
\label{2.11}
\end{equation} 
where $t={\omega}\E^{\I \varphi}$, $\lambda{\to}+\infty$. Relations  \eqref{2.8}, \eqref{2.10} and \eqref{2.11} imply \eqref{2.5}
for the case $\lambda{\to}\infty$.
In the same way  we deduce \eqref{2.5} for the case $\lambda{\to}-\infty$. 

Relation \eqref{2.6} follows from \eqref{2.5} and from the equality
(see \cite{Sakh3})
\begin{equation}
-\I \lambda\widetilde{s_1}(\lambda)=-\I \lambda\widetilde{s}(\lambda)+\E^{\I \lambda\omega}
\big( -\ln {\omega}\big)^{-\beta}+o(1),\quad\lambda{\to}\infty.
\label{2.12}
\end{equation}
The theorem is proved.\end{proof}
Thus, the function $-\I \lambda\widetilde{s}(\lambda)$ is bounded. Hence, we obtain the following assertion (see \cite{Sakh3}).
\begin{Cy}\label{Corollary 2.2} 
If $\Re{\beta}{\geq}0$ and $0<\omega<1$, then
 the operator $S_{\beta}$  defined by formula \eqref{2.2} is bounded in the space
 $L^{2}(0,\omega)$.\end{Cy}
 Further we need the following asymptotic  relation (see \cite{Sakh3}): 
\begin{equation}
\big\| S_{\beta}\E^{-\I x\lambda}+\I {\lambda}\widetilde{s_1}(\lambda)\E^{-\I x\lambda}\big\|
{\to}0,\quad \lambda{\to}\infty.
\label{2.13}
\end{equation}
Using \eqref{2.6} and \eqref{2.13} we derive two statements below.
\begin{Cy}\label{Corollary 2.3}
If $\Re{\beta}=0$ and $|z|=1$, then the points $z$ belong
to the spectrum of the operator $S_{\beta}$.
\end{Cy}
\begin{Cy}\label{Corollary 2.4}
If $\Re{\beta}=0$,\, $\nu>0$,\, $\gamma=\beta+\nu$, then in the space $L^{2}(0,\omega)$ we have
\begin{equation}
S_{\gamma}{\to}S_{\beta},\quad \nu{\to}0 \quad \mathrm{(strong \,\, convergence)}.
\label{2.14}\end{equation}
\end{Cy}

\paragraph{2.}
Taking into account \eqref{1.3} we can reformulate Corollaries \ref{Corollary 2.2}--\ref{Corollary 2.4} for the operator $V_{\beta}$.

\begin{Cy}\label{Corollary 2.5} 
If $\Re{\beta}{\geq}0$ and $0<\omega<1$, then
 the operator $V_{\beta}$   is bounded in the space
 $L^{2}(0,\infty)$.
\end{Cy}
  
\begin{Cy}\label{Corollary 2.6}
If $\Re{\beta}=0$ and $|z|=1$, then the points $z$ belong
to the spectrum of the operator $V_{\beta}$.
\end{Cy}

\begin{Cy}\label{Corollary 2.7}
If $\Re{\beta}=0$, $\nu>0$, $\gamma=\beta+\nu$, then in the space $L^{2}(0,\infty)$ we have
\begin{equation}
V_{\gamma}{\to}V_{\beta},\quad \nu{\to}0, \quad \mathrm{(strong \,\, convergence)}.
\label{2.15}\end{equation}
\end{Cy}

\paragraph{3.}
 Let us introduce the operator
\begin{equation}
T_{\beta}f=\int_{0}^{x}f(t)\big| \ln (x-t)\big|^{-\beta}\d t, \quad \beta>0,
\label{2.16}
\end{equation}
where $f(x){\in}L^{2}(0,\omega)$ and $0<\omega<1.$

\begin{Pn}\label{Proposition 2.8}
 The operator $T_{\beta}$ is such that
$T_{\beta}{\in}\sigma_{2}$ and $T_{\beta}{\notin}\sigma_{1}$.
\end{Pn}

\begin{proof}
It is obvious that $T_{\beta}{\in}\sigma_{2}$. 
In order to prove that $T_{\beta}{\notin}\sigma_{1}$ we write the operator $T_{\beta}$ in the form
\begin{equation}T_{\beta}f=\frac{\d}{\d x}\int_{0}^{x}f(t)R(x-t)\d t,\label{2.17}\end{equation}
where
\begin{equation}R(x)=\int_{0}^{x}\ln^{-\beta}(u)\d u.\label{2.18}\end{equation}
Integrating by parts the right-hand side of the expression
\begin{equation}\widetilde{s}_{1}(\lambda)=\int_{0}^{\omega}\E^{\I t\lambda}R(t)
(1-t/\omega)\d t,
\label{2.19}\end{equation}
we obtain
\begin{equation}\widetilde{s}_{1}(\lambda)=\frac{-1}{\I \lambda}\int_{0}^{\omega}\E^{\I t\lambda}\big|\ln (t)\big|^{-\beta}\d t.
\label{2.20}\end{equation}
Comparing relations \eqref{2.3} and \eqref{2.20}, we see that
${-\I \lambda}\widetilde{s}_1(\lambda )$ for the operator $T_{\beta}$ is equal to $\widetilde{s}(\lambda )$ for the operator $S_{\beta}$. Relations \eqref{2.5} and \eqref{2.20} imply that
\begin{equation}\widetilde{s}_{1}(\lambda)=\frac{\E^{\I \omega\lambda}}{\lambda^{2}}
(-\ln {\omega})^{-\beta}\big( 1+o(1)\big) ,
\quad \lambda{\to}\infty.
\label{2.21}
\end{equation}
It follows from \eqref{2.13} and \eqref{2.21} that
\begin{equation}
\sum_{n=1}^{\infty}\big| (T_{\beta}\varphi_n,\varphi_n)\big| =\infty ,
\label{2.22}
\end{equation}
where $\varphi_n(x)=\E^{\I xn/\omega}$.

Next, we use the following result \cite{FK}: \textit{
 if an operator $A$ is compact in the Hilbert space and $\varphi_j$ is an orthonormal system, then
\begin{equation}\sum_{1}^{n}\big| (A\varphi_j,\varphi_j)\big|{\leq}\sum_{j=1}^{n}s_n,
\label{2.23}
\end{equation} 
where $s_j(A)$ is the  non-decreasing  sequence of the eigenvalues of the operator} $(AA^{*})^{1/2}$.
The formulated result implies the following statement: 
\textit{if the equality}
\begin{equation}\sum_{1}^{\infty}\big| (A\varphi_j,\varphi_j)\big| =\infty\label{2.24}
\end{equation}
\textit{is valid, then} $\sum_{j=1}^{\infty}s_j=\infty$. In other words, \eqref{2.22} yields
\begin{equation}\sum_{j=1}^{\infty}s_j=\infty.
\label{2.25}
\end{equation}
Hence, the proposition is proved.
\end{proof}

\section{Friedrichs model}
\label{sec3}
Let us consider in the Hilbert space $L^{2}(0,\omega)$ the operator
\begin{equation}A_{\alpha,\beta}f=V_{\beta}QV_{\alpha}f,\quad Qf=xf,\quad \Re{\alpha}>0,
\quad \Re{\beta}>0.\label{3.1}\end{equation}
Using relations \eqref{1.1} and \eqref{1.2} we obtain
\begin{equation}A_{\alpha,\beta}f=\int_{0}^{x}f(u)U_{\alpha,\beta}(x,u)\d u,\label{3.2}\end{equation}
where
\begin{equation}U_{\alpha,\beta}(x,u)=\frac{1}{\Gamma(\alpha)\Gamma(\beta)}
\int_{0}^{x-u}E_{\beta}(x-u-y)(u+y)E_{\alpha}(y)\d y.\label{3.3}\end{equation}
We need the relation (see \cite{BaErd}):
\begin{equation}\int_{0}^{x}(x-y)^{\alpha-1}y^{\beta-1}\d y=
\frac{\Gamma(\alpha)\Gamma(\beta)}{\Gamma(\alpha+\beta)}x^{\alpha+\beta-1},\,\,
\Re{\alpha}>0,\,\,\Re{\beta}>0.
\label{3.4}
\end{equation}
Relations \eqref{1.2} and \eqref{3.4} imply that
\begin{equation}\frac{1}{\Gamma(\alpha)\Gamma(\beta)}
\int_{0}^{x-u}E_{\beta}(x-u-y)uE_{\alpha}(y)\d y=
\frac{1}{\Gamma(\alpha+\beta)}E_{\alpha+\beta}(x-u)u,\label{3.5}\end{equation}
\begin{equation}\frac{1}{\Gamma(\alpha)\Gamma(\beta)}
\int_{0}^{x-u}E_{\beta}(x-u-y)yE_{\alpha}(y)\d y=
\frac{(x-u)\alpha}{\Gamma(\alpha+\beta+1)}E_{\alpha+\beta}(x-u).\label{3.6}\end{equation}
Relations \eqref{3.3} and \eqref{3.5}, \eqref{3.6}  imply the following proposition.

\begin{Pn}\label{Proposition 3.1}
The operator $A_{\alpha,\beta}$ is defined by
\eqref{3.2}, where
\begin{equation}U_{\alpha,\beta}(x,u)=\frac{1}{\Gamma(\alpha+\beta)}E_{\alpha+\beta}(x-u)u+
\frac{(x-u)\alpha}{\Gamma(\alpha+\beta+1)}E_{\alpha+\beta}(x-u).\label{3.7}\end{equation}
\end{Pn}Using again relations \eqref{1.2} and \eqref{3.4} we have
\begin{equation}\frac{1}{\Gamma(m)}V_{\beta}E_{m}(y)=\frac{1}{\Gamma(m+\beta)}E_{m+\beta}(x),
\quad \Re{\beta}>0,\quad \Re{m}>0.\label{3.8}\end{equation}
Hence, the equality
\begin{equation}\lim_{\beta{\to}0}V_{\beta}\left(\frac{1}{\Gamma(m)}E_{m}(y)\right) =
\frac{1}{\Gamma(m)}E_{m}(x)\label{3.9}\end{equation}
is valid.
\begin{Rk}\label{Remark 3.1}It is easy to see, that the functions $E_{m}(x)$ form a complete system in the Hilbert space $L^{2}(0,\omega)$,\end{Rk}
Thus, we proved the following statement.

\begin{Pn}\label{Proposition 3.2} 
The operator   $V_{\beta}$ strongly converges  to the identity operator $I$ when $\beta{\to}+0$.
\end{Pn}
In view of \eqref{3.1}, \eqref{3.7} and Proposition \ref{Proposition 3.2} the following theorem is valid.
\begin{Tm}\label{Theorem 3.4}
Let
$\Re{\alpha}=0$, $\alpha{\ne}0$, $\beta=-\alpha.$ Then the operator 
\begin{equation}
A_{\alpha}=V_{\alpha}^{-1}QV_{\alpha} \label{3.10}\end{equation} has the form
\begin{equation}A_{\alpha}f=xf(x)+{\alpha}\int_{0}^{x}(x-y)E_{0}(x-y)f(y)\d y,\quad
f(x){\in}L^{2}(0,\omega),\label{3.11}\end{equation}
where
\begin{equation}E_{0}(x)=\int_{0}^{\infty}\frac{1}{\Gamma(s)}\E^{-Cs}s^{-1}x^{s-1}\d s.
\label{3.12}\end{equation}\end{Tm}
It follows from \eqref{3.12} that
\begin{equation}xE_{0}(x)=\frac{1}{|\ln {x}|}\big( 1+o(1)\big) ,\quad x{\to}0.
\label{3.13}\end{equation}
Let us formulate the analogue of the Theorem \ref{Theorem 3.4} (\cite[Ch. 3, Section 4]{Sakh4}).

\begin{Tm}\label{Theorem 3.5}
Let
$\Re{\alpha}=0$, $\alpha{\ne}0.$ Then the operator \begin{equation}B_{\alpha}=J_{\alpha}^{-1}QJ_{\alpha} \label{3.14}\end{equation} has the form
\begin{equation}B_{\alpha}f=xf(x)+{\alpha}\int_{0}^{x}f(y)\d y,\quad
f(x){\in}L^{2}(0,\omega).\label{3.15}\end{equation}\end{Tm}
We note, that the operators $A_{\alpha}$ and $B_{\alpha}$ are  partial cases of the Friedrichs model \cite{LFad}.

\section{Generalized wave operators}
\label{sec4}
\paragraph{1.} 
Let us introduce the notion  of the generalized wave operators (\cite{Sakh8}).

\begin{Dn}\label{Definition 4.1}
Let the operators $A$ and $A_{0}$ act in the Hilbert space $H$, where the operator $A_{0}$ is  self-adjoint with absolutely continuous spectrum. We assume that there exists a unitary operator function $W_0(t)$
satisfying the following conditions:
\begin{enumerate}
\item 
The limits in the sense of strong convergence
\begin{equation}
W_{\pm}(A,A_{0}) =\lim_{t{\to}\pm\infty}\big(\E^{iAt}\E^{-iA_{0}t}W_{0}(t)\big)
\label{4.1}
\end{equation}
exist.
\item 
\begin{equation}
\lim_{t{\to}\pm\infty}W_{0}^{-1}(t+\tau)W_{0}(t)=I.
\label{4.2}
\end{equation}
\item 
The commutations relations  hold for arbitrary values $t$ and $\tau$:
\begin{equation}
W_{0}(t)A_{0}=A_{0}W_{0}(t),\quad W_{0}(t)W_{0}(t+\tau)=W_{0}(t+\tau)W_{0}(t).
\label{4.3}
\end{equation}
\end{enumerate}
The operators $W_{\pm}(A,A_{0})$ are named  generalized wave operators. \\
If
$W_{0}(t)=I$ then the operators $W_{\pm}(A,A_{0})$ are usual wave operators.
\end{Dn}

The formulated notions  of wave operators and generalized wave operators
are  correct and useful not only for  self-adjoint operators $A$, but for non-self-adjoint operators $A$ too.

\begin{Dn}\label{Definition 4.2}
The spectrum of non-self-adjoint operator $A$
is  absolutely continuous  if the operator $A$ can be represented in the form
\begin{equation}
A=V^{-1}A_{0}V, 
\label{4.4}
\end{equation}
where the operators $V$ and $V^{-1}$ are bounded and the operator $A_{0}$ is  self-adjoint with absolutely continuous spectrum.
\end{Dn}

\begin{Pn}\label{Proposition 4.3} 
Let conditions \eqref{4.1}--\eqref{4.3} be fulfilled.
Then
\begin{equation} W_{\pm}(A,A_{0})\E^{\I A_{0}t}=\E^{\I At} W_{\pm}(A,A_{0})
\label{4.5}\end{equation}
\end{Pn}
\begin{proof}
 As in the of self-adjoint case we use the relation
\begin{equation}W_{\pm}(A,A_{0})=\lim_{t{\to}\pm\infty}\big( \E^{\I A(t+s)}\E^{-\I A_{0}(t+s)}W_{0}(t+s)\big) .
\label{4.6}\end{equation}
 The second relation of \eqref{4.3} and \eqref{4.6} imply
 \begin{equation}W_{\pm}(A,A_{0})=\E^{\I As}W_{\pm}(A,A_{0})\E^{-\I A_{0}s}.
\label{4.7}\end{equation}
This proves the proposition.
\end{proof}

\begin{Ee}\label{Example 4.4} 
Let us consider the operator $A_{\alpha}$, where
$\Re{\alpha}=0$, $\alpha{\ne}0$, and the operator $A_{0}=Q$ (see Theorem \ref{Theorem 3.4}).
\end{Ee}
According to equality \eqref{3.10} the operator $A_{\alpha}$ has absolutely continuous spectrum.
The following statement is valid.

\begin{Pn}\label{Proposition 4.5} 
We assume that
\begin{equation} 
W_{0}(t)=\big(\ln^{-\alpha}{|t|}\big) I
\label{4.8}.
\end{equation}
We  have
\begin{equation}
W_{\pm}(A_{\alpha},Q)=V_{\alpha}^{-1}.\label{4.9}\end{equation}
\end{Pn}
\begin{proof}
Relations \eqref{4.2} and \eqref{4.3} are fulfilled. Now, we write the equality
\begin{equation}
\E^{\I A_{\alpha}t}\E^{-\I Qt}W_{0}(t)\E^{\I x\tau}=
V_{\alpha}^{-1}\E^{\I Qt}V_{\alpha}\E^{-\I Qt}W_{0}(t)\E^{\I x\tau} .
\label{4.10}
\end{equation}
Using \eqref{2.6}, \eqref{2.13}, \eqref{4.2} we obtain \eqref{4.9}. The proposition is proved.
\end{proof}

\begin{Ee}\label{Example 4.6} 
Let us consider the operator $B_{\alpha}$, where
$\Re{\alpha}=0$, $\alpha{\ne}0$, and the operator $A_{0}=Q$ (see Theorem \ref{Theorem 3.5}).
\end{Ee}
According to equality \eqref{3.14} the operator $B_{\alpha}$ has absolutely continuous spectrum.
The following statement is valid.
\begin{Pn}\label{Proposition 4.7} 
Assume that
\begin{equation} W_{0}(t)=\big( |t|^{-\alpha}\big) I\label{4.11}.\end{equation}
Then we  have
\begin{equation}
W_{\pm}(B_{\alpha},Q)=\E^{\pm(\I \alpha\pi/2)}J^{-\alpha}.\label{4.12}\end{equation}
\end{Pn}
\begin{proof}
Relations \eqref{4.2} and \eqref{4.3} are fulfilled. For the operator  $J^{\alpha}$ ($\Re{\alpha}=0$, $\alpha{\ne}0$) the following relation
\begin{equation}
-\I \lambda\widetilde{s_1}(\lambda)=|\lambda|^{-\alpha}\big( 1+o(1)\big) \E^{\pm(\I \alpha\pi/2)}
,\quad \lambda{\to}\pm\infty
\label{4.13}
\end{equation} 
holds (see \cite[formula (22)]{Sakh3}).
Now, we write the equality
\begin{equation}\E^{\I B_{\alpha}t}\E^{-\I Qt}W_{0}(t)\E^{\I x\tau}=
J^{-\alpha}\E^{\I Qt}J^{\alpha}\E^{-\I Qt}W_{0}(t)\E^{\I x\tau}.
\label{4.14}\end{equation}
Hence, using  \eqref{2.13}, \eqref{4.13} and  \eqref{4.2} we obtain \eqref{4.12}. The proposition is proved.
\end{proof}


\begin{thebibliography}{19}

\bibitem{BaErd}
H. Bateman  and A. Erdelyi, \textit{ Higher transcendental 
functions, I}, New York-London, Mc Graw-Hill Book Company, 1953.

\bibitem{LFad}
 L.D. Faddeev, \textit{On a model of Friedrichs in the theory of perturbations of the continuous spectrum,} Amer.
Math. Soc. Transl. (2) \textbf{62}, 177--203, 1967.

\bibitem{FK}
Ky Fan, \textit{Maximum properties and inequalities for the eigenvalues of completely
continuous operators}, Proc. Nat. Acad. Sci. USA \textbf{37}:11 (1951), 760--766.


\bibitem{Gokr}
I. Gohberg and M.G. Krein, \textit{Theory and applications  of  Volterra operators in Hilbert space}, Amer. Math. Soc., Providence, 1970.

\bibitem{HiPhi}
E. Hille and R. Phillips, \textit{ Functional analysis  and  semigroups}, Amer. Math. Soc., Providence, 1957.

\bibitem{Mark}
A.I. Markushevich, \textit{Theory of functions of a complex variable, 3}, AMS Chelsea Publishing,
1977  (translated from Russian).

\bibitem{Sakh8}
L.A. Sakhnovich, \textit{Dissipative operators with absolutely 
continuous spectrum}, Trans. Mosc. Math. Soc. \textbf{19} (1968),  233--297. 


\bibitem{Sakh3}
L.A. Sakhnovich, \textit{Triangular integro-differential 
operators with difference kernels}, Sib. Math. Journ. \textbf{19}:4 (1978),  871--877 (Russian).


\bibitem{Sakh4}
L.A. Sakhnovich,  \textit {Integral equations with difference
kernels on finite intervals,} 1st edition,   Operator Theory Adv. Appl., vol. {84}, Birkh\"auser, Basel--Boston--Berlin, 1996.

 \bibitem{Sakh5}
L.A. Sakhnovich, \textit{On the notion and important classes of (S+N)-triangular operators},
 arXiv:1501.02831.

\bibitem{Tit}
E.C. Titchmarsh, \textit{Introduction to the theory of Fourier  integrals}, Oxford,~1937.


\end{thebibliography}
\end{document}